\documentclass[11pt]{article}

\usepackage{amsmath, amssymb, amscd, amsthm}

\newtheorem{theorem}{Theorem}

\def\T{\mbox{$\mathcal{I}$}}
\def\Mod{\mbox{Mod}}
\def\Z{\mbox{$\mathbb{Z}$}}
\def\R{\mbox{$\mathbb{R}$}}
\def\SU{\mbox{SU}}
\def\HH{\mbox{H}}

\title{The Casson invariant and the word metric on the Torelli group}

\author{Nathan Broaddus, Benson Farb, Andrew Putman\thanks{The first
two authors are supported in part by the NSF}}

\date{July 13, 2007}

\begin{document}

\maketitle

\begin{abstract}
We bound the value of the Casson invariant of any integral homology
$3$-sphere $M$ by a constant times the distance-squared to the
identity, measured in any word metric on the Torelli group $\T$,
of the element of $\T$ associated to any Heegaard splitting of
$M$.  We construct examples which show this bound is asymptotically
sharp.
\end{abstract}

\section{Introduction}
\label{section:introduction}

\noindent
The {\em Casson invariant} $\lambda(M)\in \Z$ is a fundamental and
well-studied invariant of integral homology $3$-spheres $M$.  
Roughly speaking, $\lambda(M)$ is half the algebraic number of conjugacy classes of irreducible
representations of $\pi_1(M)$ into $\SU (2)$.  See \cite{AM} for a
thorough exposition of the Casson invariant.  

The {\em mapping class group} $\Mod_g$ of a closed, orientable, genus $g$ surface
$\Sigma_g$ is the group of homotopy classes of
orientation-preserving homeomorphisms of $\Sigma_g$.  The subgroup of
$\Mod_g$ consisting of elements acting trivially on $\HH_1(\Sigma_g;\Z)$
is called the {\em Torelli group}, and is denoted by $\T_g$.  

Let $M$ be an integral homology
$3$-sphere, and let $f:\Sigma_g \to M$ be a Heegaard embedding.  For any
$\phi \in \T_g$, denote by $M_\phi$ the homology $3$-sphere obtained by
cutting $M$ along $f(\Sigma_g)$ and gluing back the resulting two
handlebodies $M^+$ and $M^-$ along their boundaries 
via the homeomorphism $\phi$.  Note that any integral homology
$3$-sphere can be obtained from $M=S^3$ in this way.  

In this note we give a sharp asymptotic bound on $|\lambda(M_\phi)|$ in
terms of the word metric on $\T_g$.  To explain our result, we fix $g >
2$ and pick once and for all a finite set $S$ of generators for $\T_g$;
the fact that $\T_g$ is finitely generated when $g>2$ is a deep result
of D. Johnson (see
\cite{Jo2}).  Denote by 
$\Vert \cdot \Vert$ the induced word norm on $\T_g$; i.e.\ $\Vert \phi
\Vert$ is the length of the shortest word in $S^{\pm 1}$ which equals
$\phi$.  Different choices of finite generating sets for $\T_g$ give
word norms whose ratios are bounded by a constant. For a fixed Heegaard embedding $f:\Sigma_g \to M$,  Morita \cite{Mo} has defined a kind of 
{\em normalized Casson invariant}
$\lambda_f:\T_g\to \Z$ via 
$$\lambda_f(\phi):=\lambda(M_\phi)-\lambda(M).$$
In particular, if $M=S^3$ and $h:\Sigma_g \to S^3$ is the unique genus $g$ Heegaard embedding then $\lambda(S^3)=0$, so the normalized Casson invariant $\lambda_h$ satisfies $\lambda_h(\phi) = \lambda(S^3_\phi)$.  

\medskip
\begin{theorem}
\label{theorem:main}
Let $M$ be an oriented integral homology $3$-sphere, let $g>2$, and let
$f:\Sigma_g\to M$ be a Heegaard embedding.  Then there exists a constant $C>0$ so that
$|\lambda_f(\phi)| \leq C \Vert \phi \Vert^2$
for every $\phi \in \T_g$.  This bound is sharp in the sense that there exists an
infinite set $\{\phi_n  \} \subset \T_g$ and a constant $K>0$ so that $|\lambda_f(\phi_n)| \geq K \Vert \phi_n \Vert^2$ for all $n$.
\end{theorem}

\medskip
\noindent
For the case $g=2$, the Torelli group $\T_2$ is not finitely generated \cite{MM}.

\section{Morita's formula}
Our proof of Theorem \ref{theorem:main} relies
in an essential way on a beautiful formula due to Morita \cite{Mo} for $\lambda_f(\phi)$, 
which we now explain (following \S 4 of \cite{Mo}).  This formula
measures the extent to which
$\lambda_f$ fails to be a homomorphism.  This failure is encoded as a function 
$\delta_f:\T_g\times\T_g\to \Z$ defined as follows.  Let $\T_{g,1}$
denote the Torelli group of an oriented, genus $g$ surface with one
boundary component $\Sigma_{g,1}$.  In other words, $\T_{g,1}$ is the group of
homotopy classes of orientation-preserving homeomorphisms of
$\Sigma_{g,1}$ which fix the boundary pointwise, modulo
homotopies which do the same and where the homeomorphisms act trivially on
$H:=\HH_1(\Sigma_g;\Z)$.  Gluing a disc to $\partial \Sigma_{g,1}$ induces a natural surjective
homomorphism $\pi: \T_{g,1}\to \T_g$, and there is a 
corresponding commutative diagram of {\em Johnson homomorphisms} 
(c.f.\ \cite{Jo1} for discussions of these homomorphisms $\tau$ and their remarkable properties):
$$\begin{CD}
\T_{g,1} @>{\tau}>> \wedge^3 H\\
@V{\pi}VV          @VVV\\
\T_{g}   @>{\tau}>> \wedge^3 H / H
\end{CD}$$
The map $f:\Sigma_g\to M$ induces 
homomorphisms $H\to \HH_1(M^\pm;\Z)$
whose kernels we denote by $H^+$ and $H^-$, respectively.  It is then
easy to see that $H^+\otimes\R$ and $H^-\otimes\R$ are maximal isotropic
subspaces of the symplectic vector space $H\otimes\R$, and that 
$$H=H^+\oplus H^-.$$
Moreover, since $M$ is an integral homology 3-sphere, there is a symplectic basis 
$$\{x_1,\ldots,x_g,y_1,\ldots ,y_g\}$$
for $H$ with $x_i\in H^+$
and $y_i\in H^-$.  Now, given any two $\phi,\psi\in\T_g$, choose
any lifts $\tilde{\phi}, \tilde{\psi}$ to $\T_{g,1}$.  Using the
obvious basis for $\wedge^3 H$ coming from our choice
of basis for $H$, we can write
\begin{align*}
\tau(\tilde{\phi}) &= \Big[\sum_{i<j<k}a_{ijk} \ y_i\wedge y_j \wedge y_k \, \Big] + \mbox{other terms}, \\
\tau(\tilde{\psi}) &= \Big[\sum_{i<j<k}b_{ijk} \ x_i\wedge x_j \wedge x_k \, \Big] + \mbox{other terms}
\end{align*}
for some $a_{ijk},b_{ijk}\in\Z$.  Morita defines
$$\delta_f(\phi,\psi)=\sum_{i<j<k}a_{ijk}b_{ijk}$$
and proves that $\delta_f(\phi,\psi)$ does not depend on either the
choice of lifts $\tilde{\phi}, \tilde{\psi}$ or the choice of symplectic
basis for $H$.  Morita then proves, as Theorem~4.3 of \cite{Mo}, 
that the following formula holds for all $\phi,\psi\in\T_g$:

\begin{equation}
\label{eq:morita:formula}
\lambda_f(\phi\psi)=\lambda_f(\phi)+\lambda_f(\psi)+2\delta_f(\phi,\psi).
\end{equation}

\section{Proof of Theorem \ref{theorem:main}}
Let $\{x_1, \ldots, x_g, y_1, \ldots, y_g\}$ be the standard basis for $H := \HH_1(\Sigma_g;\Z)$
discussed in the previous section.  For any vector $v \in \wedge^3 H$,
we denote by $\ell(v)$ the maximum of the absolute values of
the coefficients of $v$ with respect to the induced basis for $\wedge^3 H$.

We want to relate $\lambda_f(\phi)$ to the word length of $\phi$ in $\T_g$, but
Morita's formula (\ref{eq:morita:formula}) is computed using elements of
$\T_{g,1}$, not of $\T_g$.  To address this point, we first recall that 
gluing a disk to $\partial\Sigma_{g,1}$ induces an exact sequence 
$$
\begin{CD}
1 @>>> \pi_1(T^1\Sigma_g) @>>> \T_{g,1} @>{\pi}>> \T_g @>>> 1,
\end{CD}
$$
where $T^1\Sigma_g$ is the unit tangent bundle of $\Sigma_g$.  For each
generator $s \in S$ of $\T_g$, choose a single lift $\tilde{s}\in\T_{g,1}$, and denote by 
$\widetilde{S}$ the union of these elements.  We can then choose as
a generating set for $\T_{g,1}$ the set $\widetilde{S}$ together
with a finite generating set for $\pi_1(T^1\Sigma_g)$.  With these choices
of generating sets, we note that each $\phi\in\T_g$ has some lift
$\tilde{\phi}$ so that
\begin{equation}
\label{eq:liftnorm}
\Vert\tilde{\phi}\Vert_{\T_{g,1}}=\Vert\phi\Vert_{\T_g}.
\end{equation}
This equality follows by writing out $\phi$ as a product of elements of $S$, 
then lifting generator by generator.  Henceforth whenever we choose
a lift of an element $\phi\in\T_g$, we will always choose a lift
$\tilde{\phi}$ satisfying (\ref{eq:liftnorm}).  The main point is that
in computing with (\ref{eq:morita:formula}), we are allowed to choose
any lifts, since Morita proves that $\delta_f(\phi,\psi)$ does not depend
on the choice of lifts. Thus we can choose lifts which do not alter word
length.  

Now since $\widetilde{S}$ is finite, there exists $C_1$ so that 
\begin{equation}
\label{eq:mform1}
\ell(\tau(\tilde{s}))\leq C_1 \ \ \ 
\mbox{for all $s\in \widetilde{S}^{\pm 1}$}.
\end{equation}
Since $\tau$ is a homomorphism to the abelian group $\wedge^3H$, it follows from
(\ref{eq:mform1}) that
\begin{equation}
\label{eq:mform2}
\ell(\tau(\tilde{\phi}))\leq C_1\Vert\tilde{\phi}\Vert \ \ \ \mbox{for all 
$\tilde{\phi}\in\T_{g,1}$}.
\end{equation}
Finally, consider $\phi,\psi\in\T_g$ together with lifts $\tilde{\phi},\tilde{\psi}$ 
satisfying (\ref{eq:liftnorm}).  If $a_{ijk}$
(resp. $b_{ijk}$) are the coordinates of $\tau(\tilde{\phi})$
(resp. $\tau(\tilde{\psi})$) as in the previous section, then
\begin{equation}
\label{eq:mform3}
\begin{split}
|\delta_f(\phi,\psi)| &=    \Big|\!\sum_{i<j<k}a_{ijk}b_{ijk} \, \Big| 
                               \leq \Big|\!\sum_{i<j<k}\ell(\tau(\tilde{\phi})) \ell(\tau(\tilde{\psi}))\, \Big| \\
                      &\leq \sum_{i<j<k} C_1^2  \ \Vert\phi\Vert \ \Vert\psi\Vert
                               \leq C_2 \Vert\phi\Vert\ \Vert\psi\Vert
\end{split}
\end{equation}
\noindent
where $C_2={\textstyle {2g \choose 3}}C_1^2$.  

Now given any $\phi\in\T_g$, write 
$\phi=s_1\cdots s_n$, where each $s_i$ is an element of $S^{\pm 1}$ and where
$n=\Vert\phi\Vert$.  An iterated use of Morita's formula 
(\ref{eq:morita:formula}) gives 
\begin{equation}
\label{eq:keyform}
\begin{split}
\lambda_f(\phi) &= \lambda_f(s_1)+\lambda_f(s_2\cdots s_n)+2\delta_f(s_1,s_2\cdots s_n)\\
                &= \lambda_f(s_1)+\lambda_f(s_2)+ \lambda_f(s_3\cdots s_n)+2\delta_f(s_1,s_2\cdots s_n) +2\delta_f(s_2,s_3\cdots s_n)\\
                &\ \ \ \ \ \vdots\\
                &=\sum_{i=1}^n\lambda_f(s_n) + 2\sum_{i=1}^{n-1}\delta_f(s_i,s_{i+1} \cdots s_n).
\end{split}
\end{equation}
Since $S$ is finite, there exists $C_3>0$ so that
$|\lambda_f(s)|\leq C_3$ for every $s\in S$.  For some $C>0$, we thus have
\begin{equation*}
\begin{split}
|\lambda_f(\phi)| &\leq \sum_{i=1}^n|\lambda_f(s_n)| + 2\sum_{i=1}^{n-1}|\delta_f(s_i,s_{i+1} \cdots s_n)| \\
                  &\leq C_3n + 2\sum_{i=1}^{n-1}C_2 \cdot 1\cdot (n-i) \\
                  &\leq Cn^2 = C \Vert\phi\Vert^2.
\end{split}
\end{equation*}
The first claim of the theorem follows.  

We now consider the second claim.  Johnson proved (see, e.g.\ \cite{Jo1})
that the homomorphisms $\tau$ are surjective.  Hence there exists 
some $\nu \in \T_g$ so that for some lift $\tilde{\nu}\in\T_{g,1}$ we have
$$\tau(\tilde{\nu})= x_1\wedge x_2\wedge x_3 + y_1\wedge y_2\wedge y_3,$$
and hence 
\begin{equation}
\label{eq:mform5}
\tau(\tilde{\nu}^n)= n (x_1\wedge x_2\wedge x_3) + n (y_1\wedge y_2\wedge
y_3).
\end{equation}
Note that the choice of $\nu$ depends in a nontrivial way 
on the Heegaard embedding $f:\Sigma_g\to M$, so $\nu$ is not given
explicitly.  By equation (\ref{eq:keyform}), we have
\begin{equation}
\label{eq:mform6}
\lambda_f(\nu^n)=\sum_{i=1}^n\lambda_f(\nu) + 2\sum_{i=1}^{n-1}\delta_f(\nu,\nu^{n-i}).
\end{equation}
Now let $K_1=|\lambda_f(\nu)|$, which is a constant since $\nu$ is fixed.  
By (\ref{eq:mform5}) and the definition of $\delta_f$, we have for any 
$m>0$ that $\delta_f(\nu,\nu^m)=m$.  Thus by equation (\ref{eq:mform6}) there is some $N$ such that for all $n \geq N$ we have 
\begin{align*}
|\lambda_f(\nu^n)| &=    \Big|\sum_{i=1}^n \lambda_f(\nu) + 2\sum_{i=1}^{n-1}(n-i)\,\Big| \\
                   &\geq 2\sum_{i=1}^{n-1}(n-i) - \sum_{i=1}^n K_1 \geq K_2n^2
\end{align*}
for some $K_2>0$.  If $\Vert \nu \Vert=K_3$, then clearly $\Vert \nu^n \Vert \leq K_3n$.  Thus
$$|\lambda_f(\nu^n)| \geq  K_2n^2  \geq  \frac{\displaystyle
K_2}{\displaystyle K_3^2} \Vert \nu^n \Vert^2 \ \ \ \ \ \ \mbox{for all $n \geq N$}.$$
Setting $K=\frac{\displaystyle
K_2}{\displaystyle K_3^2}$ we get the desired infinite set $\{ \nu^n | n \geq N \} \subset \T_g$ establishing the asymptotic tightness of the upper bound.

\noindent
Dept. of Mathematics, University of Chicago\\
5734 S. University Ave.\\
Chicago, IL 60637\\

\noindent
E-mails: {\tt broaddus@math.uchicago.edu} (Nathan Broaddus),\\
{\tt farb@math.uchicago.edu} (Benson Farb),\\
{\tt andyp@math.uchicago.edu} (Andrew Putman)

\end{document}